\documentclass[a4paper,12pt,final]{amsart}
\usepackage{times,a4wide,mathrsfs,amssymb,amsmath,amsthm,enumerate,xypic,tikzsymbols,dsfont}

\newcommand{\C}{\mathbb{C}}

\newcommand{\QQ}{\mathbb{Q}}
\newcommand{\NN}{\mathbb{N}}
\newcommand{\PP}{\mathbb{P}}

\newcommand{\OO}{\mathcal O}

\newcommand{\XX}{\mathcal X}
\newcommand{\YY}{\mathcal Y}

\newcommand{\Zz}{\mathcal Z}

\newcommand{\MM}{\mathcal M}
\newcommand{\NNN}{\mathcal N}

\newcommand{\wt}{\widetilde}

\newcommand{\one}{\mathds{1}}

\DeclareMathOperator{\ima}{Im}

\newtheorem{theorem}{Theorem}[section]
\newtheorem{claim}[theorem]{Claim}
\newtheorem{lemma}[theorem]{Lemma}

\newtheorem{corollary}[theorem]{Corollary}
\newtheorem{proposition}[theorem]{Proposition}

\newtheorem{convention}{Conventions}

\newtheorem{nonumbering}{Theorem}

\newtheorem{nonumberingc}{Corollary}

\theoremstyle{definition}
\newtheorem{remark}[theorem]{Remark}
\newtheorem{definition}[theorem]{Definition}

\newtheorem{notation}[theorem]{Notation}

\newtheorem{nonumberingt}{Acknowledgments}

\begin{document}

\author[Robert Laterveer]
{Robert Laterveer}

\address{Institut de Recherche Math\'ematique Avanc\'ee,
CNRS -- Universit\'e 
de Strasbourg,\
7 Rue Ren\'e Des\-car\-tes, 67084 Strasbourg CEDEX,
FRANCE.}
\email{robert.laterveer@math.unistra.fr}

\title{Algebraic cycles and intersections of a quadric and a cubic}

\begin{abstract} 
Let $Y$ be a smooth complete intersection of a quadric and a cubic in $\PP^n$, with $n$ even. We show that $Y$ has a multiplicative Chow--K\"unneth decomposition, in the sense of Shen--Vial. As a 
consequence, the Chow ring of (powers of) $Y$ displays K3-like behaviour. As a by-product of the argument, we also establish a multiplicative Chow--K\"unneth decomposition for the resolution of singularities of a general nodal cubic hypersurface of even dimension.
\end{abstract}


\thanks{Supported by ANR grant ANR-20-CE40-0023.}

\maketitle

\section{Introduction}

Given a smooth projective variety $Y$ over $\C$, let $A^i(Y):=CH^i(Y)_{\QQ}$ denote the Chow groups of $Y$ (i.e. the groups of codimension $i$ algebraic cycles on $Y$ with $\QQ$-coefficients, modulo rational equivalence \cite{F}). The intersection product defines a ring structure on $A^\ast(Y)=\bigoplus_i A^i(Y)$, the {\em Chow ring\/} of $Y$.
In the case of K3 surfaces, this ring structure has a remarkable property:

\begin{theorem}[Beauville--Voisin \cite{BV}]\label{K3} Let $S$ be a K3 surface. 
The $\QQ$-subalgebra
  \[  R^\ast(S):=  \bigl\langle  A^1(S), c_j(S) \bigr\rangle\ \ \ \subset\ A^\ast(S) \]
  injects into cohomology under the cycle class map.
  \end{theorem}

Motivated by the particular behaviour of Chow rings of K3 surfaces and of abelian varieties \cite{Beau}, Beauville \cite{Beau3} has conjectured that for certain special varieties, the Chow ring should admit a multiplicative splitting. To make concrete sense of Beauville's elusive ``splitting property'' conjecture, Shen--Vial \cite{SV} introduced the concept of {\em multiplicative Chow--K\"unneth decomposition\/}. It is something of a challenge to understand the class of special varieties admitting such a decomposition: abelian varieties, K3 surfaces, cubic hypersurfaces, hyperelliptic curves and some Fano varieties are in this class (but not all Fano varieties), cf. subsection \ref{ss:mck} below.

The main result of the present paper aims to contribute to this program:

\begin{nonumbering}[=Theorem \ref{main}] Let $Y\subset\PP^{n}(\C)$ be a smooth complete intersection of a quadric and a cubic, where $n$ is even. Then $Y$ has a (generically defined) multiplicative Chow--K\"unneth decomposition.
\end{nonumbering}

%
%
%

The assumption on $n$ is there for good reason: indeed, the statement of the theorem is {\em false\/} for $n=3$ (cf. Remark \ref{!}). The question what happens for $n>3$ odd is an interesting open problem.

In proving Theorem \ref{main}, we rely on the connection between a general $Y$ and a nodal cubic $n$-fold, and we apply instances of the {\em Franchetta property\/} (cf. Section \ref{s:fr} below). 

Theorem \ref{main} implies that the Chow ring of $Y$ behaves just like that of K3 surfaces, i.e. the image of the intersection product is as small as possible:

\begin{nonumberingc}[=Corollary \ref{cor}] Let $Y\subset\PP^{n}$ be as in Theorem \ref{main}. Then
  \[  \ima\Bigl( A^i(Y)\otimes A^j(Y)\to A^{i+j}(Y)\Bigr) =\QQ[h^{i+j}]\ \ \ \forall\ i,j>0\ \]
  (where $h\in A^1(Y)$ denotes the hyperplane class).
\end{nonumberingc}

Using Theorem \ref{main}, we can also prove a result concerning the {\em tautological ring\/}, which is a certain subring of the Chow ring of powers of $Y$:

\begin{nonumberingc}[=Corollary \ref{cor1}] Let $Y\subset\PP^{n}$ be as in Theorem \ref{main}, and $m\in\NN$. Let
  \[ R^\ast(Y^m):=\Bigl\langle (p_i)^\ast(h), (p_{ij})^\ast(\Delta_Y)\Bigr\rangle\ \subset\ \ \ A^\ast(Y^m)   \]
  be the $\QQ$-subalgebra generated by pullbacks of the polarization $h\in A^1(Y)$ and pullbacks of the diagonal $\Delta_Y\in A^{n-2}(Y\times Y)$. 
  The cycle class map induces injections
   \[ R^\ast(Y^m)\ \hookrightarrow\ H^\ast(Y^m,\QQ)\ \ \ \ \forall\ m\le 2\dim H^{n-2}(Y,\QQ)-1  \ .\]
   \end{nonumberingc}

Again, this is similar to existing results for hyperelliptic curves and for K3 surfaces (cf. \cite{Ta2}, \cite{Ta}, \cite{Yin} and Remark \ref{tava} below).

 Theorem \ref{main} also has consequences for the Chow ring of resolutions of nodal cubic hypersurfaces:
 
 \begin{nonumberingc}[=Corollary \ref{cor3}] Let $\bar{Z}\subset\PP^{n+1}(\C)$ be a general nodal cubic hypersurface, where $n$ is even. Then there is a resolution of singularities $Z\to\bar{Z}$ such that $Z$ has a multiplicative Chow--K\"unneth decomposition. 
   \end{nonumberingc}

 \vskip0.6cm

\begin{convention} In this article, the word {\sl variety\/} will refer to a reduced irreducible scheme of finite type over $\C$. A {\sl subvariety\/} is a (possibly reducible) reduced subscheme which is equidimensional. 

{\bf All Chow groups are with rational coefficients}: we will denote by $A_j(Y)$ the Chow group of $j$-dimensional cycles on $Y$ with $\QQ$-coefficients; for $Y$ smooth of dimension $n$ the notations $A_j(Y)$ and $A^{n-j}(Y)$ are used interchangeably. 
The notations $A^j_{hom}(Y)$ and $A^j_{AJ}(Y)$ will be used to indicate the subgroup of homologically trivial (resp. Abel--Jacobi trivial) cycles.
For a morphism $f\colon X\to Y$, we will write $\Gamma_f\in A_\ast(X\times Y)$ for the graph of $f$.

The contravariant category of Chow motives (i.e., pure motives with respect to rational equivalence as in \cite{Sc}, \cite{MNP}) will be denoted 
$\MM_{\rm rat}$.
\end{convention}

\section{Preliminaries}

\subsection{Intersections of a quadric and a cubic} 

\begin{proposition}
\label{ifk} Let $\bar{Z}$ be the nodal cubic $n$-fold defined as
  \[ \bar{Z}:=  \bigl\{ y\in\PP^{n+1}\ \big\vert\   f(x_0,\ldots,x_n) + q(x_0,\ldots,x_n)\cdot x_{n+1}=0\bigr\}\ \ \ \subset\ \PP^{n+1},\]
  where $f$ and $g$ are general polynomials of degree 3 resp. 2. Let $\phi\colon Z\to\bar{Z}$ be the blow-up with center the node $P=[0:0:\ldots:0:1]$. Then
  $Z$ is isomorphic to the blow-up $\psi$ of $\PP^n$ with center
  \[ Y:= \bigl\{  y\in \PP^n\ \big\vert\ f(y)=q(y)=0\bigr\}\ \ \ \subset\ \PP^n\ .\]
  \end{proposition}
  
  \begin{proof} This is probably well-known. The case $n=7$ is contained in \cite[Section 5.7]{IFK}; the general case is similar. The point is that projection from the node $P$ defines a birational map $p$
  from $\bar{Z}$ to $\{x_{n+1}=0\}\cong\PP^n$. The blow-up $\phi$ provides a resolution of indeterminacy of this birational map; there is a commutative diagram
    \[\begin{array}[c]{ccccc}
         &&Z&&\\
         &&&&\\
       &{}^{\scriptstyle \phi  }\swarrow\ \ \ &&\ \ \ \searrow{}^{\scriptstyle \psi}&\\
       &&&&\\
       \bar{Z}&& \stackrel{p}{\dashrightarrow}&&\PP^n\\
       \end{array}\]
      \end{proof}
      
      \begin{remark} Let notation be as in Proposition \ref{ifk}.
 In case $n=4$, $Y$ is a K3 surface associated to a nodal cubic fourfold; this case is related to singular hyperk\"ahler varieties as shown in \cite{Chr}, cf. notably \cite[Lemma 3.3]{Chr} and also \cite[Chapter 6 section 1.2]{Huy}.
 
 In case $n=7$, $Y$ is a fivefold of ``CY3-type'' (i.e. the Hodge diamond of $Y$ looks like that of a Calabi--Yau threefold, cf. \cite{IFK}, \cite{IM}), and Proposition \ref{ifk} relates $Y$ to a nodal cubic sevenfold (which is another variety of CY3-type).  
   \end{remark}

 \subsection{MCK decomposition}
\label{ss:mck}

\begin{definition}[Murre \cite{Mur}] Let $X$ be a smooth projective variety of dimension $n$. We say that $X$ has a {\em CK decomposition\/} if there exists a decomposition of the diagonal
   \[ \Delta_X= \pi^0_X+ \pi^1_X+\cdots +\pi_X^{2n}\ \ \ \hbox{in}\ A^n(X\times X)\ ,\]
  such that the $\pi^i_X$ are mutually orthogonal idempotents and $(\pi_X^i)_\ast H^\ast(X,\QQ)= H^i(X,\QQ)$.
  
  (NB: ``CK decomposition'' is shorthand for ``Chow--K\"unneth decomposition''.)
\end{definition}

\begin{remark} The existence of a CK decomposition for any smooth projective variety is part of Murre's conjectures \cite{Mur}, \cite{J4}. 
\end{remark}

\begin{definition}[Shen--Vial \cite{SV}] Let $X$ be a smooth projective variety of dimension $n$. Let $\Delta_X^{sm}\in A^{2n}(X\times X\times X)$ be the class of the small diagonal
  \[ \Delta_X^{sm}:=\bigl\{ (x,x,x)\ \vert\ x\in X\bigr\}\ \subset\ X\times X\times X\ .\]
  An {\em MCK decomposition\/} is a CK decomposition $\{\pi_X^i\}$ of $X$ that is {\em multiplicative\/}, i.e. it satisfies
  \[ \pi_X^k\circ \Delta_X^{sm}\circ (\pi_X^i\times \pi_X^j)=0\ \ \ \hbox{in}\ A^{2n}(X\times X\times X)\ \ \ \hbox{for\ all\ }i+j\not=k\ .\]
  
 (NB: ``MCK decomposition'' is shorthand for ``multiplicative Chow--K\"unneth decomposition''.) 
  
  \end{definition}
  
  \begin{remark} The small diagonal (seen as a correspondence from $X\times X$ to $X$) induces the {\em multiplication morphism\/}
    \[ \Delta_X^{sm}\colon\ \  h(X)\otimes h(X)\ \to\ h(X)\ \ \ \hbox{in}\ \MM_{\rm rat}\ .\]
 Suppose $X$ has a CK decomposition
  \[ h(X)=\bigoplus_{i=0}^{2n} h^i(X)\ \ \ \hbox{in}\ \MM_{\rm rat}\ .\]
  By definition, this decomposition is multiplicative if for any $i,j$ the composition
  \[ h^i(X)\otimes h^j(X)\ \to\ h(X)\otimes h(X)\ \xrightarrow{\Delta_X^{sm}}\ h(X)\ \ \ \hbox{in}\ \MM_{\rm rat}\]
  factors through $h^{i+j}(X)$.
  
  If $X$ has an MCK decomposition, then setting
    \[ A^i_{(j)}(X):= (\pi_X^{2i-j})_\ast A^i(X) \ ,\]
    one obtains a bigraded ring structure on the Chow ring: that is, the intersection product sends $A^i_{(j)}(X)\otimes A^{i^\prime}_{(j^\prime)}(X) $ to  $A^{i+i^\prime}_{(j+j^\prime)}(X)$.
    
      It is expected that for any $X$ with an MCK decomposition, one has
    \[ A^i_{(j)}(X)\stackrel{??}{=}0\ \ \ \hbox{for}\ j<0\ ,\ \ \ A^i_{(0)}(X)\cap A^i_{hom}(X)\stackrel{??}{=}0\ ;\]
    this is related to Murre's conjectures B and D, that have been formulated for any CK decomposition \cite{Mur}.

  The property of having an MCK decomposition is severely restrictive, and is closely related to Beauville's ``splitting property conjecture'' \cite{Beau3}. 
  To give an idea: hyperelliptic curves have an MCK decomposition \cite[Example 8.16]{SV}, but the very general curve of genus $\ge 3$ does not have an MCK decomposition \cite[Example 2.3]{FLV2}. As for surfaces: a smooth quartic in $\PP^3$ has an MCK decomposition, but a very general surface of degree $ \ge 7$ in $\PP^3$ should not have an MCK decomposition \cite[Proposition 3.4]{FLV2}.
For more detailed discussion, and examples of varieties with an MCK decomposition, we refer to \cite[Section 8]{SV}, as well as \cite{V6}, \cite{SV2}, \cite{FTV}, \cite{37}, \cite{38}, \cite{39}, \cite{46}, \cite{40}, \cite{FLV2}, \cite{44}, \cite{g8}, \cite{NOY}.
   \end{remark}
   
   \begin{proposition}[Shen--Vial \cite{SV2}]\label{blowup} Let $M$ be a smooth projective variety, and let $f\colon\wt{M}\to M$ be the blow-up with center a smooth closed subvariety
$N\subset M$. Assume that
\begin{enumerate}

\item $M$ and $N$ have an MCK decomposition;

\item the Chern classes of the normal bundle $\NNN_{N/M}$ are in $A^\ast_{(0)}(N)$;

\item the graph of the inclusion morphism $N\to M$ is in $A^\ast_{(0)}(N\times M)$;

\item the Chern classes $c_j(T_M)$ are in $A^\ast_{(0)}(M)$.

\end{enumerate}
Then $\wt{M}$ has an MCK decomposition, the Chern classes $c_j(T_{\wt{M}})$ are in $A^\ast_{(0)}(\wt{M})$, and the graph $\Gamma_f$ is in $A^\ast_{(0)}(  \wt{M}\times M)$.
\end{proposition}

\begin{proof} This is \cite[Proposition 2.4]{SV2}. (NB: in loc. cit., $M$ and $N$ are required to have a {\em self-dual\/} MCK decomposition; however, the self-duality is actually a redundant hypothesis, cf. \cite[Section 6]{FV}.)
%
\end{proof}

 \section{The Franchetta property}
 \label{s:fr}
 
 \subsection{Definition}
 
 \begin{definition} Let $\XX\to B$ be a smooth projective morphism, where $\XX, B$ are smooth quasi-projective varieties. We say that $\XX\to B$ has the {\em Franchetta property in codimension $j$\/} if the following holds: for every $\Gamma\in A^j(\XX)$ such that the restriction $\Gamma\vert_{X_b}$ is homologically trivial for the very general $b\in B$, the restriction $\Gamma\vert_b$ is zero in $A^j(X_b)$ for all $b\in B$.
 
 We say that $\XX\to B$ has the {\em Franchetta property\/} if $\XX\to B$ has the Franchetta property in codimension $j$ for all $j$.
 \end{definition}
 
 This property is studied in \cite{PSY}, \cite{BL}, \cite{FLV}, \cite{FLV3}.
 
 \begin{definition} Given a family $\XX\to B$ as above, with $X:=X_b$ a fiber, we write
   \[ GDA^j_B(X):=\ima\Bigl( A^j(\XX)\to A^j(X)\Bigr) \]
   for the subgroup of {\em generically defined cycles}. 
   (In a context where it is clear to which family we are referring, the index $B$ will sometimes be suppressed from the notation.)
  \end{definition}
  
  With this notation, the Franchetta property amounts to saying that $GDA^\ast(X)$ injects into cohomology, under the cycle class map. 
 
\subsection{The families}

  \begin{notation}\label{not} Let 
   $\YY \to B$
   denote the universal family of smooth $n-2$-dimensional complete intersections of a quadric and a cubic over
   $B \subset \bar{B}:=\PP H^0\bigl(\PP^{n},\OO_{\PP^{n}}(3)\oplus \OO_{\PP^n}(2)\bigr) $, the corresponding parameter space.
    \end{notation}
  
  \begin{notation}\label{not2} Let $B$ be as in Notation \ref{not}, and let $\bar{\Zz}\subset \PP^{n+1}\times B$ denote the universal family of nodal cubic $n$-folds defined by an equation
  \[f(x_0,\ldots,x_n)+q(x_0,\ldots,x_n)\cdot x_{n+1}=0 \] as in Proposition \ref{ifk}.
  
  Let $\phi\colon \wt{\PP}^{n+1}\to\PP^{n+1}$ be the blow-up with center $P=[0:0:\ldots:0:1]$, and let $\Zz\subset \wt{\PP}^{n+1}\times B$ denote the blow-up of $\bar{\Zz}$ with center $p_1^{-1}(P)$. Let $B_0\subset B$ be the (non-empty) Zariski open such that the fibers of $\Zz\to B_0$ are smooth.
   \end{notation}
  
  \begin{lemma}\label{F} The varieties $\YY$ and $\Zz$ are smooth.
   \end{lemma}
   
   \begin{proof}   Let $\bar{\YY}\supset \YY$ denote the Zariski closure of $\YY$ in $\PP^n\times \bar{B}$.
  As $\OO_{\PP^{n}}(3)\oplus \OO_{\PP^n}(2)$ is base-point free, the projection $\YY\to\PP^n$ has the structure of a projective bundle, hence is smooth.  
   
  As for $\Zz$, according to Proposition \ref{ifk} the Zariski closure of ${\Zz}$ is isomorphic to the blow-up of $\PP^n\times \bar{B}$ with center $\bar{\YY}$, hence is smooth.
     \end{proof}

%

  \subsection{Franchetta for $Y$}
  
 \begin{proposition}\label{f1} Let $\YY\to B$ be the universal family of smooth complete intersections of bidegree $(2,3)$ in $\PP^n$ (Notation \ref{not}).
 The family $\YY\to B$ has the Franchetta property, for any $n$.
 \end{proposition}
 
 \begin{proof} This is easy. We have already seen (proof of Lemma \ref{F}) that $\bar{\YY}\to\PP^n$ is a projective bundle. Using the projective bundle formula, and reasoning as in
 \cite{PSY} and \cite{FLV}, this implies that
   \[ GDA^\ast_B(Y)=\ima\bigl( A^\ast(\PP^n)\to A^\ast(Y)\bigr)=\langle h\rangle\ ,\]
   where $h\in A^1(Y)$ is the restriction of a hyperplane section. It follows that $GDA^\ast_B(Y)$ injects into cohomology.
   \end{proof}

  \subsection{Franchetta for $Y^2$}
  
  \begin{proposition}\label{f2} Let $\YY\to B$ be as above, and assume $n$ even. Then the family $\YY^{2/B}\to B$ has the Franchetta property.
  \end{proposition}
  
  \begin{proof} We note that $\OO_{\PP^n}(2)$ and $\OO_{\PP^n}(3)$  are (not only base-point free but even) very ample, which means that the set-up verifies condition $(\ast_2)$ of \cite{FLV}. Then applying the stratified projective bundle argument as in loc. cit. (in the precise form of \cite[Proposition 2.6]{FLV3}), we find that
    \[ GDA^\ast(Y^2)=\langle h, \Delta_Y\rangle\ .\]
    Let us now check that the right-hand side injects into cohomology. Using Lemma \ref{ok} below, $GDA^j(Y^2)$ for $j\not=n-2$ is {\em decomposable\/}, i.e. 
    \[ GDA^j(Y^2)\ \ \subset\ GDA^\ast(Y)\otimes GDA^\ast(Y)\ \ \ \forall\ j\not=n-2\ .\]
    Thus for $j\not=n-2$, the injectivity reduces to Proposition \ref{f1}.
    
    For $j=n-2$, one observes that the class of the diagonal in $H^{2n-4}(Y\times Y,\QQ)$ can not be expressed as a decomposable cycle (for otherwise $Y$ would not have any transcendental cohomology, which is absurd), and so we reduce again to Proposition \ref{f1}.
    
    \begin{lemma}\label{ok} Let $Y\subset\PP^n$ be a smooth complete intersection of a cubic and a quadric, where $n$ is even. Then
    \[  \Delta_Y\cdot (p_i)^\ast(h) = {\displaystyle\sum} {1\over 6}\, (p_1)^\ast (h^k)\cdot (p_2)^\ast (h^{n-1-k})\ \ \ \hbox{in}\ A^{n-1}(Y\times Y)\ .\]
    \end{lemma}
    
    To prove the lemma, we write $Y=Y^\prime\cap Q$, where $Y^\prime$ and $Q$ are a cubic resp. a quadric. Up to shrinking the base $B$, we may assume $Y^\prime$ and $Q$ are smooth. Let $\iota\colon Y\hookrightarrow Q$ denote the inclusion morphism. Then $\Delta_Y\cdot (p_i)^\ast(h)$ is equal to 
      \[ {1\over 3}\, {}^t \Gamma_\iota\circ \Gamma_\iota = {1\over 3}\,(\iota\times\iota)^\ast (\Delta_Q)\ \ \ \hbox{in}\ A^{n-1}(Y\times Y)\ .\]
      But $Q$, being an odd-dimensional smooth quadric, has $A^j(Q)=\QQ$ for all $j$ and so $\Delta_Q$ can be written as 
      \[     \Delta_Q=\sum {1\over 2}\, (p_1)^\ast (h^k)\cdot (p_2)^\ast (h^{n-1-k})\ \ \ \hbox{in}\ A^{n-1}(Q\times Q)\ .\]
      
      The lemma, and the proposition, are proven. 
   \end{proof}
   
   \begin{remark}\label{!} The equality of Lemma \ref{ok} is remarkable, because this equality may fail for $n$ odd.
  Indeed, in case $Y\subset\PP^{m}$ is a smooth {\em hypersurface\/} (of any degree), there is equality
    \begin{equation}\label{form} \Delta_Y\cdot (p_i)^\ast(h)= {\displaystyle\sum} a_k\, (p_1)^\ast (h^k)\cdot (p_2)^\ast (h^{m-k})\ \ \ \hbox{in}\ A^{m}(Y\times Y)\ ,\end{equation}
    with $a_k\in\QQ$,
as follows from the excess intersection formula. On the other hand, in case $Y\subset\PP^m$ is a complete intersection of codimension at least 2,
in general there is {\em no equality\/} of the form \eqref{form}. Indeed, let $C$ be a very general curve of genus $g\ge 4$. The Faber--Pandharipande cycle
  \[  FP(C):= \Delta_C\cdot (p_j)^\ast(K_C) - {1\over 2g-2} K_C\times K_C\ \ \ \in A^2(C\times C)\ \ \ \ \ (j=1,2) \]
  is homologically trivial but non-zero in $A^2(C\times C)$ \cite{GG}, \cite{Yin0} (this cycle $FP(C)$ is the ``interesting 0-cycle'' in the title of \cite{GG}). In particular, for the very general complete intersection $Y\subset\PP^3$ of bidegree $(2,3)$, the cycle
  \[ FP(Y):=  \Delta_Y\cdot (p_j)^\ast(h) - {1\over 6} h \times h\ \ \ \in A^2(Y\times Y) \]
  is homologically trivial but non-zero, and so there cannot exist an equality of the form \eqref{form} for $Y$.
  
  This is intimately related to MCK decompositions. Indeed, if the curve $Y$ has an MCK decomposition which is generically defined, then $FP(Y)\in A^2_{(0)}(Y\times Y)\cong\QQ$ and so
  $FP(Y)$ would be zero.
  \end{remark}

  \begin{corollary}\label{z2} Let $\Zz\to B_0$ be as in Notation \ref{not2}, and assume $n$ even. Then the family $\Zz^{2/B_0}\to B_0$ has the Franchetta property.  
  \end{corollary}
  
  \begin{proof} The relation between $Z$ and $Y$ of Proposition \ref{ifk} is obviously generically defined (with respect to $B_0$), and so there is a generically defined isomorphism of motives
    \begin{equation}\label{motiso}  h(Z)\ \cong\ h(Y)(-1)\oplus \bigoplus_{k=0}^n \one(-k)\ \ \ \hbox{in}\ \MM_{\rm rat}\ .\end{equation}
    This induces a commutative diagram 
     \[ \begin{array}[c]{ccc}
     GDA^j_{B_0}(Z^2)&\cong&GDA^{j-2}_{B_0}(Y^2)\oplus \bigoplus GDA^\ast_{B_0}(Y)\oplus \QQ^r\\
       &&\\
       \downarrow&&\downarrow\\
       &&\\
       H^{2j}(Z^2,\QQ)&\cong& H^{2j-4}(Y^2,\QQ)\oplus \bigoplus H^\ast(Y,\QQ)\oplus \QQ^r\\
       \end{array}
            \]
     The corollary now follows from Propositions \ref{f1} and \ref{f2}.
       \end{proof}

  \subsection{Franchetta for $Y^3$ in codimension $\le 2n-4$}
  
  \begin{proposition}\label{f3} Let $\YY\to B$ be as above, and assume $n$ even. Then the family $\YY^{3/B}\to B$ has the Franchetta property in codimension $\le 2n-4$.
  \end{proposition}
  
  \begin{proof} We exploit the relation between $Y$ and $Z$ of Proposition \ref{ifk}. The isomorphism of motives \eqref{motiso} induces a commutative diagram 
    \[  \begin{array}[c]{ccc}  GDA^j_{B_0}(Y^3)& \hookrightarrow&  GDA^{j+3}_{B_0}(Z^3)\\
         &&\\
         \downarrow&&\downarrow\\
         &&\\
         H^{2j}(Y^3,\QQ) &\hookrightarrow& H^{2j+6}(Z^3,\QQ)\\
         \end{array}\]
    This diagram implies that Proposition \ref{f3} is a consequence of the following result:
    
    \begin{proposition}\label{z3} Let $\Zz\to B_0$ be as in Notation \ref{not2}, and assume $n$ even. The family $\Zz^{3/B_0}\to B_0$ has the Franchetta property in codimension $\le 2n-1$.
  \end{proposition}
  
  It remains to prove Proposition \ref{z3}. Let $\phi\colon\wt{\PP}^{n+1}\to\PP^{n+1}$ denote the blow-up with center $P$ and exceptional divisor $E$.
  The family $\Zz\to B_0$ can be interpreted as the universal family of hypersurfaces associated to the (very ample) line bundle 
    \[ L:=\phi^\ast\OO_{\PP^{n+1}}(3)\otimes \OO_{\wt{\PP}^{n+1}}(-2E)\] 
    on $\wt{\PP}^{n+1}$.
  
  Let us now convince ourselves that $(\wt{\PP}^{n+1},L)$ has property $(\ast_3)$ of \cite{FLV} (which means that 3 distinct points of $\wt{\PP}^{n+1}$ impose 3 independent conditions on the parameter space $B_0$). Three distinct points outside of $E$ impose independent conditions on sections of $L$, because 
3 distinct points impose independent conditions on cubics in $\PP^{n+1}$. Three distinct points contained in $E$ impose independent conditions on sections of $L$, because 3 distinct points impose independent conditions on quadrics in $\PP^{n}$ (note that the intersection of $Z$ with the exceptional divisor $E\cong\PP^n$ is the quadric defined by $q=0$, where $q$ is as in Proposition \ref{ifk}). Taking a combination of points in $E$ and points outside $E$, things get even easier: the points in $E$ only impose conditions on the quadric part $q$ of the equation of $\bar{Z}$ in Proposition \ref{ifk}, whereas points outside of $E$ only impose conditions on the cubic part $f$ of the equation of $\bar{Z}$. This shows that $(\wt{\PP}^{n+1},L)$ has property $(\ast_3)$.

Applying the stratified projective bundle argument of \cite{FLV} (cf. also \cite[Proposition 2.6]{FLV3} for the precise form used here), we find that
  \[ GDA^\ast_{B_0}(Z^3)=\Bigl\langle (p_i)^\ast(h), (p_j)^\ast(E), (p_{kl})^\ast(\Delta_Z)\Bigr\rangle\ ,\]
  where $h,E\in A^1(Z)$ denote the restriction of $\phi^\ast(h)$ resp. the exceptional divisor $E$ to $Z$.
         
  We note that when $j\le 2n-1$, an element in $GDA^j(Z^3)$ cannot involve a product of 2 diagonals, and so
   \[ GDA^j(Z^3)\ \ \subset\ \sum GDA^\ast(Z^2)\otimes GDA^\ast(Z)\ \ \ \forall\ j\le 2n-1\ .\]
As these summands go to different pieces of the K\"unneth decomposition in cohomology, Proposition \ref{z3} now follows from Proposition \ref{z2}.
    \end{proof}

  \section{Main result}
  
  \begin{theorem}\label{main} Let $Y\subset\PP^{n}$ be a smooth complete intersection of a quadric and a cubic, where $n$ is even. Then $Y$ has a multiplicative Chow--K\"unneth decomposition.
  \end{theorem}

 \begin{proof} 
 
 Let us first construct a CK decomposition for $Y$. Letting $h\in A^1(Y)$ denote a hyperplane section, as before we consider
    \[   \begin{split}  
                              \pi^{2j}_Y&:= {1\over 6}\, h^{n-2-j}\times h^j\ \ \ \ (j\not={n-2\over 2})\ ,\\
                                \pi^{n-2}_Y&:= \Delta_Y-\sum_{j\not={n-2\over 2}} \pi^{2j}_Y\ \ \ \ \ \ \in\ A^{n-2}(Y\times Y)\ .\\
                                \end{split}\]
 We observe that this CK decomposition is {\em generically defined\/} with respect to the family $\YY\to B$ (Notation \ref{not}), i.e. it is obtained by restriction from
 ``universal projectors'' $\pi^j_\YY\in A^{n-2}(\YY\times_B \YY)$. (This is just because $h$ and $\Delta_Y$ are generically defined.)    
 
  Writing $h^j(Y):=(Y,\pi^j_Y)\in\MM_{\rm rat}$, we have
   \[ h^{2j}(Y) \cong \one(-j)\ \ \  \hbox{in}\ \MM_{\rm rat}\ \ \ \ \ \ (j\not={n-2\over 2})\ ,\]
 i.e. the interesting part of the motive of $Y$ is contained in $h^{n-2}(Y)$.

 What we need to prove is that this CK decomposition is MCK, i.e.
      \begin{equation}\label{this} \pi_Y^k\circ \Delta_Y^{sm}\circ (\pi_Y^i\times \pi_Y^j)=0\ \ \ \hbox{in}\ A^{2n-4}(Y\times Y\times Y)\ \ \ \hbox{for\ all\ }i+j\not=k\ ,\end{equation}
      or equivalently that
       \[   h^i(Y)\otimes h^j(Y)\ \xrightarrow{\Delta_Y^{sm}}\  h(Y) \]
       coincides with 
       \[ h^i(Y)\otimes h^j(Y)\ \xrightarrow{\Delta_Y^{sm}}\ h(Y)\ \to\ h^{i+j}(Y)  \ \to\ h(Y)\ , \]   
       for all $i,j$.
       
     We observe that the cycles in \eqref{this} are generically defined, and that the vanishing \eqref{this} holds true modulo homological equivalence. As such, the required vanishing \eqref{this} follows from the Franchetta property for $\YY^{3/B}\to B$ in codimension $2n-4$ (Proposition \ref{f3}). 
               \end{proof}

  \section{Consequences}   
  
  \subsection{Chow ring of $Y$}
  
  \begin{corollary}\label{cor} Let $Y\subset\PP^{n}$ be as in Theorem \ref{main}. Then
  \[  \ima\Bigl( A^i(Y)\otimes A^j(Y)\to A^{i+j}(Y)\Bigr) =\QQ[h^{i+j}]\ \ \ \forall\ i,j>0\ .\]
\end{corollary}

\begin{proof} (NB: This argument is similar to \cite[Theorem 1.0.1]{Diaz}, which is about cubic hypersurfaces.)

We begin by noting that when $n=4$, the complete intersection $Y$ is a K3 surface and so Corollary \ref{cor} is already known from the seminal work of Beauville--Voisin \cite{BV}. We now suppose that $n\ge 6$, and so $Y$ is a Fano variety.

Let us consider the {\em modified small diagonal\/} 
  \[ \begin{split} \Gamma_3:= \Delta^{sm}_Y - {1\over 6}\Bigl( (p_{12})^\ast(\Delta_Y)\cdot (p_3)^\ast(h) +  (p_{13})^\ast(\Delta_Y)\cdot (p_2)^\ast(h)+  (p_{23})^\ast(\Delta_Y)\cdot (p_1)^\ast(h)\Bigr)&\\ +
    {\displaystyle\sum_{i+j+k=2n-4}} a_{ijk} (p_1)^\ast(h^i)\cdot (p_2)^\ast(h^j)\cdot (p_3)^\ast(h^k)\ \ \ \in\ A^{2n-4}(Y\times Y\times& Y)\ ,\\
    \end{split}\]
    where $a_{ijk}\in\QQ$. We make the following claim:
    
    \begin{claim}\label{cl} There exist $a_{ijk}\in\QQ$ such that
      \[  \Gamma_3=0\ \ \ \hbox{in}\ A^{2n-4}(Y\times Y\times Y)\]
      for all $Y$ as in Theorem \ref{main}.
      \end{claim}
      
      The claim implies the corollary, as can be seen by letting $\Gamma_3$ act on $\alpha\times\beta\in A^{i+j}(Y\times Y)$ (cf. \cite[Lemma 2.0.2]{Diaz}; this uses that $Y$ is Fano so that $A_0(Y)\cong\QQ$). To establish the claim, 
 we reason as in \cite[Proof of Proposition 2.8]{FLV3}: the MCK decomposition (plus the fact that $Y$ has transcendental cohomology, plus the relation of Lemma \ref{ok}) yields an identity of the form
      \[ \begin{split} \Delta^{sm}_Y={1\over 6}\Bigl( (p_{12})^\ast(\Delta_Y)\cdot (p_3)^\ast(h) +  (p_{13})^\ast(\Delta_Y)\cdot (p_2)^\ast(h)+  (p_{23})^\ast(\Delta_Y)\cdot (p_1)^\ast(h)\Bigr)&\\ + P\bigl( (p_1)^\ast(h),(p_2)^\ast(h), (p_3)^\ast(h)\bigr)\ \ \ \hbox{in}\ A^{2n-4}&(Y^3)\ ,\\
      \end{split}\]
   where $P$ is a symmetric polynomial with $\QQ$-coefficients.   
   
   (Alternatively, one could also establish Claim \ref{cl} by using the Franchetta property for $Y^3$ in codimension $2n-4$ (Proposition \ref{f3}); thus one is reduced to finding $a_{ijk}$ such that
   the claim is verified modulo homological equivalence; this can be done as in \cite[Lemma 2.0.1]{Diaz}.)      
     \end{proof}
     
  \begin{remark} In case $n=4$, $Y$ is a K3 surface and in this case Corollary \ref{cor} follows from the seminal work of Beauville--Voisin \cite{BV}. 
  
  For $n$ large, Corollary \ref{cor} seems to be non-trivial, because it is (expected but) not known that $A^j(Y)\cong\QQ$ for all $j<(n-2)/2$.
 \end{remark}  

%
%

 \subsection{The tautological ring}
 
 \begin{corollary}\label{cor1} Let $Y\subset\PP^{n}$ be a smooth dimensionally transverse intersection of a quadric and a cubic, where $n$ is even. Let $m\in\NN$. Let
  \[ R^\ast(Y^m):=\Bigl\langle (p_i)^\ast(h), (p_{ij})^\ast(\Delta_Y)\Bigr\rangle\ \subset\ \ \ A^\ast(Y^m)   \]
  be the $\QQ$-subalgebra generated by pullbacks of the polarization $h\in A^1(Y)$ and pullbacks of the diagonal $\Delta_Y\in A^{n-2}(Y\times Y)$. (Here $p_i$ and $p_{ij}$ denote the various projections from $Y^m$ to $Y$ resp. to $Y\times Y$).
  The cycle class map induces injections
   \[ R^\ast(Y^m)\ \hookrightarrow\ H^\ast(Y^m,\QQ)\ \ \ \hbox{for\ all\ }m\le 2b-1  \ ,\]
   where $b:=\dim H^{n-2}(Y,\QQ)$.
   
   Moreover, $R^\ast(Y^m)\to H^\ast(Y^m,\QQ)$ is injective for all $m$ if and only if $Y$ is Kimura finite-dimensional \cite{Kim}.
   \end{corollary}

\begin{proof} This is directly inspired by the analogous result for cubic hypersurfaces \cite[Section 2.3]{FLV3}, which in turn is inspired by analogous results for hyperelliptic curves \cite{Ta2}, \cite{Ta} (cf. Remark \ref{tava} below) and for K3 surfaces \cite{Yin}.

As in \cite[Section 2.3]{FLV3}, let us write $o:={1\over 6} h^{n-2}\in A^{n-2}(Y)$, and
  \[ \tau:= \Delta_Y - {1\over 6}\, \sum_{j=0}^{n-2}  h^j\times h^{n-2-j}\ \ \in\ A^{n-2}(Y\times Y) \]
  (this cycle $\tau$ is nothing but the projector on the motive $h^{n-2}_{prim}(Y)$).
Moreover, let us write 
  \[ \begin{split}   o_i&:= (p_i)^\ast(o)\ \ \in\ A^{n-2}(Y^m)\ ,\\
                        h_i&:=(p_i)^\ast(h)\ \ \in \ A^1(Y^m)\ ,\\
                         \tau_{ij}&:=(p_{ij})^\ast(\tau)\ \ \in\ A^{n-2}(Y^m)\ .\\
                         \end{split}\]
We now define the $\QQ$-subalgebra
  \[ \bar{R}^\ast(Y^m):=\Bigl\langle o_i, h_i, \tau_{ij}\Bigr\rangle\ \ \ \subset\ H^\ast(Y^m,\QQ)\ \ \ \ \ (1\le i\le m,\ 1\le i<j\le m)\ ; \]
 this is the image of $R^\ast(Y^m)$ in cohomology. One can prove (just as \cite[Lemma 2.11]{FLV3} and \cite[Lemma 2.3]{Yin}) that the $\QQ$-algebra $ \bar{R}^\ast(Y^m)$
  is isomorphic to the free graded $\QQ$-algebra generated by $o_i,h_i,\tau_{ij}$, modulo the following relations:
    \begin{equation}\label{E:X'}
			o_i\cdot o_i = 0, \quad h_i \cdot o_i = 0,  \quad 
			h_i^{n-2} =6\,o_i\,;
			\end{equation}
			\begin{equation}\label{E:X2'}
			\tau_{ij} \cdot o_i = 0 ,\quad \tau_{ij} \cdot h_i = 0, \quad \tau_{ij} \cdot \tau_{ij} = (b-1)\, o_i\cdot o_j
			\,;
			\end{equation}
			\begin{equation}\label{E:X3'}
			\tau_{ij} \cdot \tau_{ik} = \tau_{jk} \cdot o_i\,;
			\end{equation}
			\begin{equation}\label{E:X4'}
			\sum_{\sigma \in \mathfrak{S}_{b}} 
			\hbox{sign}(\sigma) 
			\prod_{i=1}^{b} \tau_{i, b+\sigma(i)} = 0\, ,
			\end{equation}
			
To prove Corollary \ref{cor1}, it suffices to check that these relations are verified modulo rational equivalence.
The relations \eqref{E:X'} take place in $R^\ast(Y)$ and so they follow from the Franchetta property for $Y$. 
The relations \eqref{E:X2'} take place in $R^\ast(Y^2)$. The first and the last relations are trivially verified, because one may assume $n>4$ and so $Y$ is Fano and one has
$A^{2n-4}(Y^2)=\QQ$. As for the second relation of \eqref{E:X2'}, this is the relation of Lemma \ref{ok}.
   
   Relation \eqref{E:X3'} takes place in $R^\ast(Y^3)$ and follows from the MCK decomposition (Theorem \ref{main}). Indeed, we have
   \[  \Delta_Y^{sm}\circ (\pi^{n-2}_Y\times\pi^{n-2}_Y)=   \pi^{2n-4}_Y\circ \Delta_Y^{sm}\circ (\pi^{n-2}_Y\times\pi^{n-2}_Y)  \ \ \ \hbox{in}\ A^{2n-4}(Y^3)\ ,\]
   which (using Lieberman's lemma) translates into
   \[ (\pi^{n-2}_Y\times \pi^{n-2}_Y\times\Delta_Y)_\ast    \Delta_Y^{sm}  =   ( \pi^{n-2}_Y\times \pi^{n-2}_Y\times\pi^{2n-4}_Y)_\ast \Delta_Y^{sm}                                                         \ \ \ \hbox{in}\ A^{2n-4}(Y^3)\ ,\]
   which means that
   \[  \tau_{13}\cdot \tau_{23}= \tau_{12}\cdot o_3\ \ \ \hbox{in}\ A^{2n-4}(Y^3)\ .\]
   
  Finally, relation \eqref{E:X4'}, which takes place in $R^\ast(Y^{2b})$
   is the Kimura finite-dimensionality relation \cite{Kim}: in the notation of loc. cit., relation \eqref{E:X4'} is written as   
    \[   \wedge^{b} \pi^{n-2, prim}_Y = 0\ \ \ \hbox{in}\ H^{b(2n-4)}(Y^{b}\times Y^{b},\QQ)\ ,\]
    which
  expresses the fact that ($\dim H^{n-2}_{prim}(Y,\QQ)=b-1$ and so)
    \[ \wedge^{b} H^{n-2}_{prim}(Y,\QQ)=0\ .\]
 Assuming that $Y$ is Kimura finite-dimensional, this relation is also verified modulo rational equivalence.
%
%
%
 \end{proof}

\begin{remark}\label{tava} Given any curve $C$ and an integer $m\in\NN$, one can define the {\em tautological ring\/}
  \[ R^\ast(C^m):=  \Bigl\langle  (p_i)^\ast(K_C),(p_{ij})^\ast(\Delta_C)\Bigr\rangle\ \ \ \subset\ A^\ast(C^m) \]
  (where $p_i, p_{ij}$ denote the various projections from $C^m$ to $C$ resp. $C\times C$).
  Tavakol has proven \cite[Corollary 6.4]{Ta} that if $C$ is a hyperelliptic curve, the cycle class map induces injections
    \[  R^\ast(C^m)\ \hookrightarrow\ H^\ast(C^m,\QQ)\ \ \ \hbox{for\ all\ }m\in\NN\ .\]
   On the other hand, there are many (non hyperelliptic) curves for which the tautological rings $R^\ast(C^2)$ or $R^\ast(C^3)$ do {\em not\/} inject into cohomology (this is related to the non-vanishing of the Faber--Pandharipande cycle of Remark \ref{!}; it is also
   related to the non-vanishing of the Ceresa cycle, cf. \cite[Remark 4.2]{Ta} and \cite[Example 2.3 and Remark 2.4]{FLV2}). 
   
Corollary \ref{cor1} shows that the tautological ring of $Y$ behaves similarly to that of hyperelliptic curves, K3 surfaces and cubic hypersurfaces. 
\end{remark}

%

\subsection{The MCK decomposition for $Z$}

\begin{corollary}\label{cor3} Let $\bar{Z}\subset\PP^{n+1}(\C)$ be a general nodal cubic, where $n$ is even. Then there is a resolution of singularities $Z\to\bar{Z}$ such that $Z$ has an MCK decomposition. 
\end{corollary}

\begin{proof} For general $\bar{Z}$, the blow-up $Z\to\bar{Z}$ with center the node is non-singular. To prove the corollary, we interpret $Z$ as the blow-up of $\PP^n$ with center $Y$ a smooth complete intersection of a cubic and a quadric (Proposition \ref{ifk}), and apply Proposition \ref{blowup} (with $M=\PP^n$ and $N=Y$). Let us check all conditions of Proposition \ref{blowup} are satisfied. In view of Theorem \ref{main}, conditions (1) and (4) are satisfied.
Conditions (2) and (3) are about generically defined cycles on $Y$, and so they follow from Proposition \ref{f1}.
\end{proof}

%
%

 \vskip1cm
\begin{nonumberingt} Thanks to Yoyo from bellenana.fr.
\end{nonumberingt}

\end{document}